\documentclass[12pt,twoside]{article}
\usepackage{graphicx}

\usepackage{amsbsy,amssymb,amsmath,amsfonts}

\usepackage{arydshln}
\allowdisplaybreaks
\usepackage{cite}
\usepackage{tikz}
\usepackage{graphics}
\def\bpsp{\begin{pspicture}}
\def\epsp{\end{pspicture}}

\newtheorem{theorem}{Theorem}[section]
\newtheorem{remark}[theorem]{Remark}
\newtheorem{example}[theorem]{Example}
\newtheorem{lemma}[theorem]{Lemma}
\newtheorem{corollary}[theorem]{Corollary}
\newtheorem{definition}[theorem]{Definition}
\newtheorem{proposition}[theorem]{Proposition}

\newtheorem{note}{Note}
\newtheorem{case}{Case}
\newtheorem{conjecture}[theorem]{Conjecture}
\newtheorem{question}[theorem]{Question}

\newcommand{\bea}{\begin{eqnarray}}
\newcommand{\eea}{\end{eqnarray}}
\newcommand{\beq}{\begin{eqnarray*}}
\newcommand{\eeq}{\end{eqnarray*}}

\def\m4{\mbox{\rm ~(mod $4$)}}

\def \bd{\begin{definition}}
\def \ed{\end{definition}}
\def \bqu{\begin{question}}
\def \equ{\end{question}}
\def \bcc{\begin{conjecture}}
\def \ecc{\end{conjecture}}
\def \bt{\begin{theorem}}
\def \et{\end{theorem}}
\def \bl{\begin{lemma}}
\def \el{\end{lemma}}
\def \bc{\begin{corollary}}
\def \ec{\end{corollary}}
\def \be{\begin{equation}}
\def \ee{\end{equation}}
\def \ben{\begin{enumerate}}
\def \een{\end{enumerate}}
\def \ba{\begin{array}}
\def \ea{\end{array}}
\def \bp{\begin{proposition}}
\def \ep{\end{proposition}}
\def \bx{\begin{example}}
\def \ex{\end{example}}
\def \br{\begin{remark}}
\def \er{\end{remark}}
\def \bdsc{\begin{description}}
\def \edsc{\end{description}}

\def \bn{\begin{case}}
\def \en{\end{case}}
\def \bnt{\begin{note}}
\def \ent{\end{note}}
\def\1{1\!\!1}

\def\mm2{\mbox{\rm ~(mod $2$)}}
\def\m4{\mbox{\rm ~(mod $4$)}}

\def\m{\mu}

\def\1{\textbf{1}}
\def\0{\textbf{0}}

\begin{document}
\title{On normalized Laplacian eigenvalues of power graphs associated to finite cyclic groups}
\author{ Bilal A. Rather$ ^{a} $, S. Pirzada$ ^{b} $,  T. A. Chishti$ ^c $, Ahmad M. Alghamdi$^{d}$  \\
$^{a,b}${\em  Department of Mathematics, University of Kashmir, Srinagar, India}\\
$^{c}${\em  Directorate of Distance Education, University of Kashmir, Srinagar, India}\\
$^{d}${\em  Department of Mathematical Sciences, Umm Alqura University, Saudi Arabia}\\
$^a$bilalahmadrr@gmail.com,~~$^b$pirzadasd@kashmiruniversity.ac.in\\
$ ^c $tachishti@uok.edu.in,~~amghamdi@uqu.edu.sa}
\date{}

\pagestyle{myheadings} \markboth{Bilal, Pirzada, Chishti, Alghamdi}{On normalized Laplacian eigenvalue of power graphs}
\maketitle
\vskip 3mm

\noindent{\footnotesize \bf Abstract.} For a simple connected graph $ G $ of order $ n $, the normalized Laplacian is a square matrix of order $ n $, defined as $\mathcal{L}(G)= D(G)^{-\frac{1}{2}}L(G)D(G)^{-\frac{1}{2}}$, where $ D(G)^{-\frac{1}{2}} $ is the diagonal matrix whose $ i$-th diagonal entry is $ \frac{1}{\sqrt{d_{i}}} $. In this article, we find the normalized Laplacian eigenvalues of the joined union of regular graphs in terms of the adjacency eigenvalues and the eigenvalues of quotient matrix associated with graph $ G $. For a finite group $\mathcal{G}$, the power graph $\mathcal{P}(\mathcal{G})$ of a group $ \mathcal{G} $ is defined as the simple graph in which two distinct vertices are joined by an edge if and only if one is the power of other. As a consequence of the joined union of graphs, we investigate the normalized Laplacian eigenvalues of power graphs of finite cyclic group $ \mathbb{Z}_{n}. $
\vskip 3mm



\section{Introduction}

A simple graph is denoted by $G(V(G),E(G))$, where $V(G)=\{v_{1},v_{2},\ldots,v_{n}\}$ is its vertex set and $E(G)$ is its edge set. The \textit{order} of $G$ is  $n=|V(G)|$ and \textit{size} is $m=|E(G)|.$ The \textit{neighborhood} of a vertex  $v$ in $ G $, denoted by $ N(v) $, is the set of all those vertices of $ G $ which are adjacent to $ v. $ The \textit{degree}  $d_{G}(v)$ (or $d_v$) of a vertex $v$ is the number of vertices in $ G $ that are incident to $ v $. The adjacency matrix, denoted by $ A(G) $, is defined by
\begin{equation*}
A(G)=\begin{cases}
1 & \text{if} ~ v_{i}\sim v_{j},\\
0 & \text{otherwise},

\end{cases}
\end{equation*}
where $v_{i} \sim v_{j} $ denotes $v_i$ is adjacent to $v_j$ in  $ G $. The eigenvalues of $ A(G) $ are denoted by $ \lambda_{i} $ and are called the adjacency eigenvalues of $ G. $ Let $D(G)={diag}(d_1, d_2, \dots, d_n)$ be the diagonal matrix of vertex degrees $d_i=d_{G}(v_i)$, $i=1,2,\dots,n$, associated to $G$. The real symmetric and positive semi-definite matrix $L(G)=D(G)-A(G)$ is the Laplacian matrix and its eigenvalues are known as Laplacian eigenvalues of  $G$.  More literature about adjacency and Laplacian matrix can be found in \cite{cds}.

The normalized Laplacian matrix of a graph $ G $, denoted by $ \mathcal{L}(G) $, is defined as
\begin{equation*}
\mathcal{L}(G)=\begin{cases}
1 & \text{if} ~ v_{i}=v_{j} ~\text{and} ~ v_{i}\neq 0,\\
\dfrac{-1}{\sqrt{d_{v_{i}}d_{v_{j}}}} & \text{if}~ v_{i}\sim v_{j},\\
0 & \text{otherwise}.

\end{cases}
\end{equation*}
This matrix was introduced by Chung \cite{chung} to study the random walks of $ G $ and is equivalently defined as $\mathcal{L}(G)= D^{-\frac{1}{2}}L(G)D^{-\frac{1}{2}}$, where $ D^{-\frac{1}{2}} $ is the diagonal matrix whose $ i$-th diagonal entry is $ \frac{1}{\sqrt{d_{i}}} $. Clearly, $ \mathcal{L}(G) $ is real symmetric and positive semi-definite matrix. Its eigenvalues are real and are known as normalized Laplacian eigenvalues. We denote normalized Laplacian eigenvalues by $ \lambda_{i}(\mathcal{L}) $ and order them as $ 0=\lambda_{1}(\mathcal{L}) \leq \lambda_{2}(\mathcal{L})\leq \dots \leq \lambda_{n}(\mathcal{L})=2$. In certain situations, normalized Laplacian is a natural tool that works better than adjacency and Laplacian matrices. More literature about $\mathcal{L}(G)$ can be seen in \cite{cavers, db, db1, sun} and the references therein.

As usual, we denote the complete graph, the bipartite graphs and the cycle graph by $ K_{n},~ K_{a,b},~ C_{n}$, respectively. For other notations and terminology, we refer to \cite{cds,sp}.

The rest of the paper is organized as follows. In Section 2, we obtain the normalized Laplacian eigenvalues of the joined union of regular graphs $G_1,G_2,\dots,G_n$ in terms of their adjacency eigenvalues and the eigenvalues of the quotient matrix associated with the joined union. In  Section 3, we discuss the normalized Laplacian eigenvalues of the power graphs of the finite cyclic groups $ \mathbb{Z}_{n} $.

\section{Normalized Laplacian eigenvalues of the joined union of graphs}

Consider the matrix
\begin{equation*}
M= \begin{pmatrix}
m_{1,1} & m_{1,2} & \cdots & m_{1,s} \\
m_{2,1} & m_{2,2} & \cdots & m_{2,s} \\
\vdots & \vdots & \ddots & \vdots \\
m_{s,1} & m_{s,2} & \cdots & m_{s,s} \\
\end{pmatrix}_{n\times n}, \end{equation*}
whose rows and columns are partitioned according to a partition $P=\{ P_{1}, P_{2},\dots , P_{m}\} $ of the set $X= \{1,2,\dots,n\}. $ The quotient matrix $ \mathcal{Q} $ of the matrix $ M $ is the $s \times s$ matrix whose entries are the average row sums of the blocks $ m_{i,j} $. The partition $ P $ is said to be \emph{ equitable} if each block $ m_{i,j} $ of $ M $ has constant row (and column) sum and in this case the matrix $ \mathcal{Q} $ is called as \emph{ equitable quotient matrix}. In general, the eigenvalues of $ \mathcal{Q}  $ interlace the eigenvalues of $ M $. In case the partition is equitable, we have following lemma.

\begin{lemma}\cite{BH, cds}
If the partition $ P $ of $ X $ of matrix $ M $ is equitable, then each eigenvalue of $ \mathcal{Q} $ is an eigenvalue of $ M. $
\end{lemma}

Let $G(V,E)$ be a graph of order $n$ and $ G_{i}(V_{i}, E_{i})$ be graphs of order $n_i,$ where $i=1,\ldots, n $. The  \textit{joined union} \cite{DS} $ G[G_{1},\ldots, G_{n}] $ is the graph $ H(W, F) $ with
\begin{equation*}
W=\bigcup_{i=1}^{n}V_{i}~~~\text{and}~~F=\bigcup_{i=1}^{n}E_{i}\bigcup\Big (\bigcup_{\{v_{i}, v_{j}\}\in E}V_{i}\times V_{j}\Big).
\end{equation*}
Equivalently, the joined union $ G[G_{1},\ldots, G_{n}] $ is obtained by joining edges from each vertex of $ G_{i} $  to every vertex of $ G_{j} $ whenever $ v_{i} \sim v_{j} $ in $ G. $ Thus, the usual join $ G_{1}\triangledown G_{2} $ is a particular case of the joined union $ K_{2}[G_1, G_2]$.

In \cite{wu}, the authors discussed the normalized Laplacian eigenvalues of $G[G_1,G_2,\dots,G_n]$ in terms of the normalized Laplacian eigenvalues of $ G_{i} $'s and the eigenvalues of another matrix using the technique of Cardosa et. al. \cite{cardosa}. Using a different approach, we will discuss the normalized Laplacian eigenvalues of $G[G_1,G_2,\dots,G_n]$ in terms of the adjacency eigenvalues of the graphs $G_1,G_2,\dots,G_n$ and the eigenvalues of the quotient matrix, where each of the $ G_{i} $ is an $ r_{i} $ regular graph.

\begin{theorem} \label{joined union}
Let $G$ be a graph of order $n$ and size $m$. Let $ G_{i}$ be $r_{i}$ regular graphs of order $ n_{i} $ having adjacency eigenvalues $\lambda_{i1}=r_{i}\geq \lambda_{i2}\geq\ldots\geq \lambda_{in_{i}}, $ where $ i=1,2, \ldots, n$.  Then the normalized Laplacian eigenvalues of the graph $ G[G_{1},\ldots, G_{n}] $ are given by
\begin{equation*}
 1-\dfrac{1}{r_{i}+\alpha_{i}}\lambda_{ik}(G_{i}), ~~for ~~ i=1,\ldots,n ~~ and ~~ k=2,3,\ldots, n_{i},
\end{equation*}
where $\alpha_i=\sum\limits_{v_j\in N_{G}(v_{i})}n_{i}$ is the sum of the orders of the graphs $G_j, j\ne i$ which correspond to the neighbours of vertex $ v_{i}\in G $. The remaining $n$ eigenvalues are given by the equitable quotient matrix $ M $ of \eqref{Qmat of joined union}.

\end{theorem}
\textbf{ Proof.}
Let $V(G)=\{ v_{1}, \ldots, v_{n}\}$ be the vertex set of $ G $ and let $V(G_i)=\{ v_{i1}, \ldots, v_{in_i}\}$ be the vertex set of the graph $G_i$, for $i=1,2,\dots,n$. Let  $H=G[G_{1}, \ldots, G_{n}]$ be the joined union  of $ r_{i} $ regular graphs $G_i,$ for $ i=1,2,\dots,n $. It is clear that the order of $ H $ is $ N=\sum\limits_{i=1}^{n}n_{i} $. Since degree of each vertex $v_{ij}\in V(H)$, is degree inside $ G_{i} $ and  the sum of orders of $ G_{j}\text{'s}, j\ne i $, which correspond to the neighbours of the vertex $ v_{i} $ in $ G $, where $1\le i\leq n$ and $1\leq j\leq n_i$, therefore, for  each $v_{ij}\in V(G_{i})$, we have
\begin{equation}\label{eq23}
d_{H}(v_{ij})=r_{i}+\sum\limits_{v_j\in N_{G}(v_i)}n_{j}=r_i+\alpha_i,
\end{equation} where $\displaystyle\alpha_i=\sum\limits_{v_j\in N_{G}(v_i)}n_{j}$. Under suitable labelling of the vertices in $ H $, the normalized Laplacian matrix of $H$ can be written as
\begin{equation*}
\mathcal{L}(H)=\begin{pmatrix}
g_1 & a(v_1,v_2)& \ldots &  a(v_1,v_n) \\
a(v_2,v_1) & g_2 & \ldots &  a(v_2,v_n)\\
\vdots &\vdots &\ddots &\vdots\\
 a(v_n,v_1) & a(v_n,v_2) &\ldots & g_n
\end{pmatrix},
\end{equation*}
where,  for $i=1,2,\ldots,n,$
\begin{equation*}
g_{i}=I_{n_{i}}-\dfrac{1}{r_{i}+\alpha_{i}}A(G_{i})~ \text{and}~ a(v_{i}, v_{j}) = \begin{cases}
\dfrac{1}{\sqrt{(\alpha_{i}+r_{i})(\alpha_{j}+r_{j})}}J_{n_{i}\times n_{j}} ,& \text{if}~  v_{i}\sim v_{j}~  \text{in}~ G\\
\textbf{0}_{n_{i}\times n_{j}}, & \text{otherwise.}
\end{cases}
\end{equation*}
$ A(G_i)$ is the adjacency matrix of $ G_{i}$, $J_{n_i\times n_{j}}$ is the matrix having all entries $ 1 $, $ \textbf{0}_{n_{i}\times n_{j}} $ is the zero matrix of order $ n_{i}\times n_{j} $ and $ I_{n_i}$ is the identity matrix of order $n_i$. \\
\indent As each $ G_{i} $ is an $ r_{i} $ regular graph, so the all one vector $e_{n_i}=(\underbrace{1,1,\dots,1}_{n_{i}})^T$  is the eigenvector of the adjacency matrix $ A(G_{i}) $ corresponding to the eigenvalue $r_{i} $ and all other eigenvectors are orthogonal to $e_{n_i}.$ Let $ \lambda_{ik}$, $2\leq k\leq n_i$, be any eigenvalue of $ A(G_{i})$ with the corresponding eigenvector $X=(x_{i1},x_{i2},\dots,x_{in_i})^T$ satisfying $e_{n_i}^TX=0.$ Clearly, the column vector $X$ can be regarded as a function defined on $ V(G_i) $ assigning the vertex $ v_{ij} $ to $ x_{ij} $, that is, $ X(v_{ij})=x_{ij} $ for $ i=1,2,\ldots,n $ and $j=1,2,\dots,n_i$. Now, consider the vector $Y=(y_{1},y_{2},\dots,y_{n})^T$, where
\begin{equation*}
y_{j}=\left\{
\begin{array}{rl}
x_{ij}&~\text{if}~ v_{ij}\in V(G_i)\\
0&~\text{otherwise.}\\
\end{array}\right.
\end{equation*}
Since, $ e_{n_{i}}^{T}X=0 $ and coordinates of the  vector $Y$ corresponding to vertices in $\cup_{j\ne i}V_j$ of $H$ are zeros,  we have
\begin{align*}
\mathcal{L}(H)Y=\begin{pmatrix}
0\\
\vdots\\
0\\
X-\dfrac{1}{r_{i}+\alpha_{i}}\lambda_{ik} X\\
0\\
\vdots\\
0
\end{pmatrix}=\left (1-\dfrac{1}{r_{i}+\alpha_{i}}\lambda_{ik}\right )Y.
\end{align*}
This shows that $Y$ is an eigenvector of $\mathcal{L}(H) $ corresponding to the eigenvalue $1-\dfrac{1}{r_{i}+\alpha_{i}}\lambda_{ik}$, for every eigenvalue $\lambda_{ik}$, $2\leq k\leq n_i$, of $A(G_i)$. In this way, we have obtained $ \sum\limits_{i=1}^{n}n_{i}-n=N-n$  eigenvalues. The remaining $ n $ normalized Laplacian eigenvalues of $ H $ are the eigenvalues of the equitable quotient matrix
\begin{equation}\label{Qmat of joined union}
M =\begin{pmatrix}
\dfrac{\alpha_{1}}{\alpha_1+r_1}& \dfrac{-n_2a_{12}}{\sqrt{(r_{1}+\alpha_{1})(r_{2}+\alpha_{2})}}&\ldots & \dfrac{-n_na_{1n}}{\sqrt{(r_{1}+\alpha_{1})(r_{n}+\alpha_{n})}}\\
\dfrac{-n_1a_{21}}{\sqrt{(r_{2}+\alpha_{2})(r_{1}+\alpha_{1})}}&\dfrac{\alpha_{2}}{\alpha_2+r_2}& \ldots & \dfrac{-n_na_{2n}}{\sqrt{(r_{2}+\alpha_{2})(r_{n}+\alpha_{n})}}\\
\vdots &\vdots  &\ddots &\vdots\\
\dfrac{-n_1a_{n1}}{\sqrt{(r_{n}+\alpha_{n})(r_{1}+\alpha_{1})}}& \dfrac{-n_2a_{n2}}{\sqrt{(r_{n}+\alpha_{n})(r_{2}+\alpha_{2})}}&\ldots &\dfrac{\alpha_{n}}{\alpha_n+r_n}
\end{pmatrix},
\end{equation}
where, for $i\ne j$,
\begin{equation*}
a_{ij}=\begin{cases}
1, & v_i\sim v_j\\
0, & \text{otherwise}.

\end{cases}
\end{equation*} $\hfill\Box$

The next observation is a consequence of Theorem \ref{joined union} and gives the normalized Laplacian eigenvalues of $K_{n_1,n_2,\dots,n_p}$.

\begin{corollary}\label{p-partite}
The normalized Laplacian eigenvalues of the complete $p$-partite graph $K_{n_1,n_2,\dots,n_p}=K_{p}[\overline{K}_{n_1},\overline{K}_{n_2},\dots,\overline{K}_{n_p}]$ with $N=\displaystyle\sum_{i=1}^{p}n_{i}$ consists of the eigenvalue $1$ with multiplicity $N-p$ and the remaining $ p $ eigenvalues are given by the matrix
\begin{equation*}
\begin{pmatrix}
1&\dfrac{-n_{2}}{\sqrt{\alpha_{1}\alpha_{2}}}&\ldots & \dfrac{-n_{p}}{\sqrt{\alpha_{1}\alpha_{p}}}\\
\dfrac{-n_{1}}{\sqrt{\alpha_{2}\alpha_{1}}}&1& \ldots & \dfrac{-n_{p}}{\sqrt{\alpha_{2}\alpha_{p}}}\\
\vdots &\vdots  &\ddots &\vdots\\
\dfrac{-n_{1}}{\sqrt{\alpha_{p}\alpha_{2}}}& \dfrac{-n_{2}}{\sqrt{\alpha_{p}\alpha_{2}}}&\ldots &1
\end{pmatrix}.
\end{equation*}
\end{corollary}
\textbf{Proof.} This follows from Theorem \ref{joined union}, by taking $ G_{i}=\overline{K}_{i} $ and $ \lambda_{ik}(G_{i})=0 $ for each $ i $ and each $ k $.  $\hfill\Box$

In particular, if partite sets are of equal size, say $ n_{1}=n_{2}=\dots= n_{p}=t $, then we have the following observation.
\begin{corollary}
Let $ G= K_{t,t,\dots,t} $ be a complete $ p$-partite graph of order $ N=pt .$ Then the normalized Laplacian eigenvalues of $ G $ consists of the eigenvalue $1 $ with multiplicity $ N-p $, the eigenvalue $\dfrac{p}{p-1}$ with multiplicity $p-1$ and the simple eigenvalue $ 0 $.
\end{corollary}
\textbf{Proof.} By Theorem \ref{joined union}, we have $ \alpha_{i}=t(p-1) $, for $ i=1,2,\dots,p $. Also, by Corollary \ref{p-partite}, we see that $ 1 $ is an eigenvalue with multiplicity $ pt-p $ and other eigenvalues are given by
\begin{equation*}
M_{p}=\begin{pmatrix}
1&\dfrac{-1}{p-1}&\ldots & \dfrac{-1}{p-1}\\
\dfrac{-1}{p-1}&1& \ldots & \dfrac{-1}{p-1}\\
\vdots &\vdots  &\ddots &\vdots\\
\dfrac{-1}{p-1}& \dfrac{-1}{p-1}&\ldots &1
\end{pmatrix}.
\end{equation*}
By simple calculations, we see that the normalized Laplacian eigenvalues of matrix $ M_{p} $ consists of the eigenvalue $ \dfrac{p}{p-1} $ with multiplicity $ p-1 $ and the simple eigenvalue $ 0. $ $\hfill\Box$

Another consequence of Theorem  \ref{joined union}, gives the normalized Laplacian eigenvalues of the join of two regular graphs.
\begin{corollary}\label{join of two graphs}
Let $ G_{i}$ be an $ r_{i} $ regular graph of order $ n_{i} $ for $i=1,2$. Let $ \lambda_{ik}, 2\leq k\leq n_{i}, i=1,2 $ be the adjacency eigenvalues of $ G_{i} $. Then the normalized Laplacian eigenvalues of $ G=G_{1}\triangledown G_{2} $ consists of the eigenvalue $ 1-\dfrac{1}{r_{1}+n_{2}}\lambda_{1k}A(G_{1}) $, $k=2,\dots, n_{1} $, the eigenvalues $  1-\dfrac{1}{r_{2}+n_{1}}\lambda_{2k}A(G_{1}), k=2,\dots, n_{2} $ and the remaining two eigenvalues are given by the quotient matrix
\begin{equation} \label{Qmat of join of two graphs}
\begin{pmatrix}
 \dfrac{n_{2}}{r_{1}+n_{2}}& \dfrac{-n_{2}}{\sqrt{(r_{1}+n_{2})(r_{2}+n_{1})}}\\
\dfrac{-n_{1}}{\sqrt{(r_{1}+n_{2})(r_{2}+n_{1})}} & \dfrac{n_{1}}{r_{2}+n_{1}}
\end{pmatrix}.
\end{equation}
\end{corollary}

Since $ G_{1} $ and $ G_{2} $ are regular graphs, we observe that the two eigenvalues of matrix \eqref{Qmat of join of two graphs} are the largest and the smallest normalized Laplacian eigenvalue of $ G=G_{1}\triangledown G_{2} $.
\begin{proposition}
The largest and the smallest normalized Laplacian eigenvalues of $ G_{1}\triangledown G_{2} $ are the eigenvalues of the matrix \eqref{Qmat of join of two graphs}.

\end{proposition}
\begin{proposition}\label{spec of join of two graphs}
\begin{itemize}
\item[ \bf{(i)} ] The normalized Laplacian eigenvalues of the complete bipartite graph $ K_{a,b}=K_{a}\triangledown K_{b} $ are
\[
\left  \{ 0,1^{[a+b-2]},2 \right \}.
 \]
\item[ \bf{(ii)} ] The normalized Laplacian eigenvalues of the complete split graph $CS_{\omega,n-\omega}=K_{\omega}\triangledown \overline{K}_{n-\omega}$, with clique number $ \omega $ and independence number $ n-\omega $ are given by
\[
\left  \{ 0,\left (\dfrac{n}{n-1}\right)^{[\omega-1]},\dfrac{2n-\omega+1}{n-1} \right \}.
 \]
\item[ \bf{(iii)}] The normalized Laplacian eigenvalues of the cone graph $ C_{a,b}=C_{a}\triangledown \overline{K}_{b} $ consists of the eigenvalues $ 1-\dfrac{1}{2+b}2\cos\left( \frac{2\pi k}{n}\right) $, where $ k=2,\dots,n-1 $, the simple eigenvalues $ 0 $ and $ \dfrac{2b+2}{b+2} $.
\item[ \bf{(iv)} ] The normalized Laplacian eigenvalues of the wheel graph $ W_{n}=C_{n-1}\triangledown K_{1} $ consists of the eigenvalues $ 1-\frac{1}{3}2\cos\left( \frac{2\pi k}{m}\right) $, where $ k=2,\dots,n-2 $, and the simple eigenvalues $ \left \{0,\dfrac{4}{3}\right \}. $
\end{itemize}
\end{proposition}
\textbf{Proof.} $ \textbf{(i)} $. This follows from Corollary \ref{join of two graphs},  by taking $ n_{1}=a,n_{2}=b, r_{1}=r_{2}=0 $ and $ \lambda_{1k}=0,$ for $ k=2,\dots,a$ and $ \lambda_{2k} =0$ for each $ k=2,\dots,b$.\\
$ \textbf{(ii)} $. We recall that the adjacency spectrum of $ K_{\omega} $ is $ \{ \omega-1,-1^{(\omega-1)} \} $. Now, the result follows from Corollary \ref{join of two graphs} by taking $ n_{1}=\omega, n_{2}=n-\omega, r_{1}=\omega-1, r_{2}=0, \lambda_{1k}=-1,$ for $ k=2,\dots,\omega $ and $ \lambda_{2k}=0$ for $k=2,3,\dots,n-\omega .$ \\
$ \textbf{(iii)} $. Since adjacency spectrum of $ C_{n} $ is $ \left\lbrace 2\cos\left( \frac{2\pi k}{n}\right) : k=1,2,\dots,n\right\rbrace  $, by taking $ n_{1}=a,n_{2}=b, r_{1}=2, r_{2}=0 $ and $ \lambda_{1k} =2\cos\left( \frac{2\pi i}{m}\right )$ for $ k=2,3,\dots,a-1 $ and $ \lambda_{2,k} $ for $ k=2,\dots,b-1 $ in Corollary \ref{join of two graphs}, we get the required eigenvalues.\\
$ \textbf{(iv)} $. This is a special case of part (iii) with $ a=n-1 $ and $ b=1 .$ $\hfill\Box$

A \emph{friendship} graph $ F_{n} $ is a graph of order $ 2n+1 $, obtained by joining $ K_{1} $ with $ n $ copies of $ K_{2} $, that is, $ F_{n}=K_{1}\triangledown (nK_{2}) =K_{1,n}[K_{1},\underbrace{K_{2},K_{2},\dots,K_{2}}_{n}]$, where $K_1$ corresponds to the root vertex (vertex of degree greater than one) in $K_{1,n}$.  In particular, replacing some of $ K_{2} $'s by $ K_{1} $'s in $F_n$ we get a \emph{firefly} type graph, denoted by $ F_{p,n-p}$ and written as
\begin{equation*}
F_{p,n-p}=K_{1,n}[K_{1},\underbrace{K_{1},K_{1},\dots,K_{1}}_{p}\underbrace{K_{2},K_{2},\dots,K_{2}}_{n-p}].
\end{equation*}
A \emph{generalized or multi-step wheel network} $ W_{a,b} $ is a graph derived from $ a $ copies of $ C_{b} $ and $ K_{1} $, in such a way that all the vertices of each $ C_{b} $ are adjacent to $ K_{1} $. Its order is $ ab+1 $ and can be written as $ W_{a,b}=K_{1}\triangledown (aC_{b})= K_{1,a}[K_{1},\underbrace{C_{b},\dots,C_{b}}_{a} ]$.

The normalized Laplacian eigenvalues of the \emph{friendship} graph $ F_{n} $, the \emph{firefly} type graph $ F_{p,n-p}$ and $ W_{a,b} $ are given by the following.
\begin{proposition}\begin{itemize}
\item[ \bf{(i)}] The normalized Laplacian eigenvalues of $ F_{n} $ are
\[ \bigg\{0, \left ( \dfrac{1}{2} \right )^{n-1}, \left ( \dfrac{3}{2} \right )^{n+1} \bigg\}. \]

\item[ \bf{(ii)}] The normalized Laplacian eigenvalues of $ F_{p,n-p} $  are
\[
\bigg\{0, \left ( \dfrac{1}{2} \right )^{[n-p-1]}, 1^{[p-1]}, \left ( \dfrac{3}{2} \right )^{[n-p]}, \dfrac{5\sqrt{2n-p}\pm\sqrt{2n+7p}}{4\sqrt{2n-p}} \bigg\}.
 \]

\item[ \bf{(iii)}] The normalized Laplacian eigenvalues of $ W_{a,b} $ consists of the eigenvalues $ 0$, the eigenvalue $ \dfrac{4}{3} $ and the eigenvalues $ 1-\dfrac{2}{3}\cos\left (\dfrac{2\pi k}{b}\right ),$ for  $ k=2,\dots,b. $

\end{itemize}
\end{proposition}
\textbf{Proof.} \textbf{(i)}. By Theorem \ref{joined union} and the definition of $ F_{n} $, we have
\[ \alpha_{1}=2n, \alpha_{2}=\dots=\alpha_{n+1}=1~~ \text{and}~~ r_{1}=0, r_{2}=\dots =r_{n+1}=1. \]
So, by Theorem \ref{joined union}, we see that $ \dfrac{3}{2} $ is the normalized Laplacian eigenvalues of $ F_{n} $ with multiplicity $ n. $ The remaining eigenvalues are given by the block matrix
 \begin{equation}\label{Qmat of Fn}
\left (\begin{array}{c|c c c c}
1 & \dfrac{-1}{\sqrt{n}}& \dots & \dfrac{-1}{\sqrt{n}}&\dfrac{-1}{\sqrt{n}}\\
\hline
\dfrac{-1}{2\sqrt{n}} & \dfrac{1}{2} & \dots & 0& 0\\
\vdots & \vdots &\ddots& \vdots & \vdots\\
\dfrac{-1}{2\sqrt{n}}& 0 & \dots &  \dfrac{1}{2} & 0\\
\dfrac{-1}{2\sqrt{n}}& 0 & \dots & 0 & \dfrac{1}{2}
\end{array}\right ).
\end{equation}
Clearly, $ \dfrac{1}{2} $ is the normalized laplacian eigenvalue of \eqref{Qmat of Fn} with multiplicity $ n-1 $ and the remaining two eigenvalues of block matrix \eqref{Qmat of Fn} are given by the quotient matrix
\[ \begin{pmatrix}
1 & \frac{-n}{\sqrt{n}}\\
\dfrac{-1}{2\sqrt{n}} & \dfrac{1}{2}

\end{pmatrix}.
 \]
\textbf{(ii)}. Since $ \alpha_{1}=p+2(n-p)=2n-p $ and $ \alpha_{2}=\dots=\alpha_{2n+1-p}=1 $, so by Theorem \ref{joined union}, with $ r_{1}=\dots=r_{p+1}=0, ~ r_{p+2}=\dots=r_{2n+1-p}=1,$ we see that $ \dfrac{3}{2} $ is the normalized Laplacian eigenvalue of $ F_{p,n-p} $ with multiplicity $ n-p $. The other normalized Laplacian eigenvalues of $ F_{p,n-p} $ are given by the block matrix
\begin{equation}\label{Qmat Fp,n-p}
\left(
\begin{array}{c| c c c| c c c}
1 & \dfrac{-1}{\sqrt{2n-p}} & \dots & \dfrac{-1}{\sqrt{2n-p}}& \dfrac{-2}{\sqrt{2(2n-p)}}& \dots & \dfrac{-2}{\sqrt{2(2n-p)}}\\
\hline
\dfrac{-1}{\sqrt{2n-p}} & 1 & \dots & 0 & 0 &\dots & 0\\
\vdots & \vdots & \ddots & \vdots & \vdots &  \dots & \vdots\\
\dfrac{-1}{\sqrt{2n-p}} & 0 & \dots & 1 & 0 & \dots & 0\\

\hline
\dfrac{-2}{\sqrt{2(2n-p)}} & 0 & \dots & 0 & \dfrac{1}{2} & \dots & 0\\
\vdots & \vdots & \dots & \vdots & \vdots &  \ddots & \vdots\\
\dfrac{-2}{\sqrt{2(2n-p)}} & 0 & \dots & 0 & 0 & \dots & \dfrac{1}{2}\\
\end{array}
\right).
\end{equation}
By simple calculations, $ 1 $ and $ \dfrac{1}{2} $ are the normalized Laplacian eigenvalues of \eqref{Qmat Fp,n-p} and the remaining eigenvalues of block matrix \eqref{Qmat Fp,n-p} are given by the quotient matrix
\begin{equation}\label{qqmat}
\begin{pmatrix}
1 & \dfrac{-p}{\sqrt{2n-p}} & \dfrac{-2(n-p)}{\sqrt{2(2n-p)}}\\
\dfrac{-1}{\sqrt{2n-p}} & 1 & 0\\
\dfrac{-1}{\sqrt{2(2n-p)}} & 0 & \dfrac{1}{2}
\end{pmatrix}.
\end{equation}
Now, it is easy to see that $ 0 $ and $ \dfrac{5\sqrt{2n-p}\pm\sqrt{2n+7p}}{4\sqrt{2n-p}} $ are the normalized Laplacian eigenvalues of quotient matrix \eqref{qqmat}.\\
\textbf{(iii)}. As in part (iii) of Proposition \ref{spec of join of two graphs}, we see that $ 1-\dfrac{2}{3}\cos\left (\dfrac{2\pi k}{b}\right ),$ for  $ k=2,\dots,b. $ are the normalized Laplacian eigenvalues of $ W_{a,b} $. The other eigenvalues are given by the block matrix
\begin{equation*}
\left( \begin{array}{c| c c c c}
1 & \dfrac{-b}{\sqrt{3ab}}& \dots & \dfrac{-b}{\sqrt{3ab}}&\dfrac{-b}{\sqrt{3ab}}\\
\hline
\dfrac{-1}{\sqrt{3ab}} & \dfrac{1}{3} & \dots & 0 &  0\\
\vdots & \vdots &\ddots& \vdots & \vdots\\
\dfrac{-1}{\sqrt{3ab}} & 0 & \dots& \dfrac{1}{3} & 0\\
\dfrac{-1}{\sqrt{3ab}}& 0 & \dots & 0 & \dfrac{1}{3}
\end{array}\right).
\end{equation*}
Now, as in part (i), $ \bigg\{0,\left (\dfrac{1}{3}\right )^{a-1}, \dfrac{4}{3}\bigg\} $ are the remaining normalized Laplacian eigenvalues of $ W_{a,b} $. $\hfill\Box$

\section{Normalized Laplacian eigenvalues of the power graphs of cyclic group $ \mathbb{Z}_{n} $ }
In this section, we consider the power graphs of finite cyclic group $ \mathbb{Z}_{n} $. As an application to Theorem \ref{joined union} and its consequences obtained in Section 2, we determine the normalized Laplacian eigenvalues of power graph of $ \mathbb{Z}_{n}. $

All groups are assumed to be finite and every cyclic group of order $ n $ is taken as isomorphic copy of integral additive modulo group $ \mathbb{Z}_{n} $ with identity denoted by $ 0 $. Let $ \mathcal{G} $ be a finite group of order $ n $ with identity  $ e $. The power graph of group $ \mathcal{G} $, denoted by $ \mathcal{P}(\mathcal{G}) $, is the simple graph with vertex set as the elements of group $ \mathcal{G} $ and two distinct vertices $ x,y\in \mathcal{G} $ are adjacent if and only if one is the positive power of the other, that is, $ x^{i}=y $ or $ y^{j}=x $, for positive integers $i,j$ with  $ 2\leq i,j\leq n $. These graphs were introduced in \cite{kelarev}, see also \cite{sen}. Such graphs have valuable applications and are related to automata theory \cite{kelarev9}, besides being useful in characterizing finite groups.  We let $ U_{n}^{*}=\{x\in\mathbb{Z}_{n}: (x,n)=1 \}\cup\{0\} $, where $ (x,n) $ denotes greatest common divisor of $ x $ and $ n $. Our other group theory notations are standard and can be taken from \cite{roman}. More work on power graphs can be seen in \cite{cameron1, sen, mehreen, survey,tamiza} and the references therein. \\

The adjacency spectrum, the Laplacian and the signless Laplacian spectrum of power graphs of finite cyclic and dihedral groups have been investigated in  \cite{banerjee,sriparna,asma, mehreen1, panda}. The normalized Laplacian eigenvalues of power graphs of certain finite groups were studied in \cite{asma1}.\\

Let $ n $ be a positive integer and $ d $ divides $ n $, written as $ d|n $. The divisor $ d $ is the proper divisor of $ n $, if $ 1<d<n . $ Let $\mathbb{G}_{n}$ be a simple graph with vertex set as the proper divisor set $ \{d_{i}: 1,n\neq d_{i}|n,  ~ 1 \leq i \leq t\} $ and edge set $\{ d_{i}d_{j}: d_{i}|d_{j}, ~ 1 \leq i< j\leq t\} $, for $ 1\leq i<j\leq t $. If the \emph{canonical decomposition} of $n$ is $ n=p_{1}^{n_{1}}p_{2}^{n_{2}}\dots p_{r}^{n_{r}} $, where $ r,n_{1},n_{2},\dots,n_{r} $ are positive integers and $ p_{1},p_{2},\dots,p_{r} $ are distinct prime numbers, then the number of divisors of $ n $ are $\prod\limits_{i=1}^{r}(n_{i}+1) $. So the order of graph  $\mathbb{G}_{n}$ is  $ |V( \mathbb{G}_{n})|=\prod\limits_{i=1}^{r}(n_{i}+1)-2. $ Also, $\mathbb{G}_{n}$ is a connected graph \cite{bilal}, provided $ n $ is neither a prime power nor the product of two distinct primes. In \cite{mehreen}, $\mathbb{G}_{n}$ is used as the underlying graph for studying the power graph of finite cyclic group $ \mathbb{Z}_{n} $ and it has been shown that for each proper divisor $ d_{i} $ of $ n $, $ \mathcal{P}(\mathbb{Z}_{n})$ has a complete subgraph of order $ \phi(d_{i}) $.  \\

The following theorem shows that $ \dfrac{n}{n-1} $ is always the normalized Laplacian eigenvalue of the power graph $ \mathcal{P}(\mathbb{Z}_{n})$.

\begin{theorem}\label{mul of alpha n/n-1 of Zn}
Let $ \mathbb{Z}_{n} $ be a finite cyclic group of order $ n\geq 3 $. Then $ \dfrac{n}{n-1} $ is normalized Laplacian eigenvalue of $ \mathcal{P}(\mathbb{Z}_{n}) $ with multiplicity at least $ \phi(n)$.
\end{theorem}
\textbf{Proof.} Let $ \mathbb{Z}_{n} $ be the cyclic group of order $ n\geq 3 $. Then the identity $ 0 $ and invertible elements of  the group $ \mathbb{Z}_{n} $ in the power graph $ \mathcal{P}(\mathbb{Z}_{n}) $ are adjacent to every other vertex in $ \mathcal{P}(\mathbb{Z}_{n}) $. Since it is well known that the number of invertible elements of $ \mathbb{Z}_{n} $ are $ \phi(n) $ in number, so the induced power graph $ \mathcal{P}(U_{n}^{*}) $ is the complete graph $ K_{\phi(n)+1} $. Thus, by Theorem \ref{mehreen}, we see that $ \mathcal{P}(\mathbb{Z}_{n})=K_{\phi(n)+1}\triangledown \mathcal{P}(\mathbb{Z}_{n}\setminus U_{n}^{*})$. By applying Corollary \ref{join of two graphs}, we get \begin{equation*}
1-\dfrac{1}{r_{1}+\alpha_{1}}(-1)=1+\dfrac{1}{\phi(n)+n-\phi(n)-1}=\dfrac{n}{n-1}
\end{equation*}
as the normalized Laplacian eigenvalue with multiplicity at least $ \phi(n) $, since $ \dfrac{n}{n-1} $ can also be the normalized Laplacian eigenvalue of quotient matrix \eqref{Qmat of join of two graphs}.$\hfill\Box$

If $ n=p^{z} $, where $ p $ is prime and $ z $ is a positive integer, then we have following observation.

\begin{corollary}\label{NL spectra of p^z}
If $ n=p^{z} $, where $ p $ is prime and $ z $ is a positive integer, then the normalized Laplacian eigenvalues of  $\mathcal{P}(\mathbb{Z}_n) $ are $ \left \{ \left (\dfrac{n}{n-1}\right )^{[(n-1)]}, ~0 \right \} $.
\end{corollary}
\textbf{Proof.} If $ n=p^z $, where $ p $ is prime and $ z $ is a positive integer, then as shown in \cite{sen}, $ \mathcal{P}(\mathbb{Z}_n)$ is isomorphic to the complete graph $ K_{n}$ and hence the result follows. $\hfill\Box$

The next observation gives the normalized Laplacian eigenvalues of $ \mathcal{P}(\mathbb{Z}_n) $, when $ n $ is the product of two primes.

\begin{corollary}\label{NL spectra of pq}
Let $ n=pq$ be the product of two distinct primes. Then the normalized Laplacian eigenvalues of $ \mathcal{P}(\mathbb{Z}_n) $ are $ \left \{ 0,\left (\dfrac{n}{n-1}\right )^{[\phi(n)]}, \left (1+\dfrac{1}{q\phi(p)}\right )^{[\phi(p)-1]}, \left (1+\dfrac{1}{p\phi(q)}\right )^{[\phi(q)-1]} \right  \} $ and the zeros of polynomial
\begin{align*}
p(x)=x\Bigg(x^2&-x\left (\dfrac{\phi(n)+1}{q\phi(p)}+\dfrac{\phi(p)+\phi(q)}{q\phi(p)+\phi(q)}+\dfrac{\phi(n)+1}{p\phi(q)}\right )+\dfrac{(\phi(n)+1)\phi(p)}{p\phi(q)(q\phi(p)+\phi(q))} \\
&+\dfrac{(\phi(n)+1)^{2}}{n\phi(n)}+\dfrac{(\phi(n)+1)\phi(q)}{q\phi(p)(q\phi(p)+\phi(q))} \Bigg).
\end{align*}

\end{corollary}
\textbf{Proof.} If $ n=pq $, where $p$ and $q$, ($ p<q) $) are primes, then $ \mathcal{P}(\mathbb{Z}_{n}) $ \cite{tamiza} can be written as
\[  \mathcal{P}(\mathbb{Z}_n)=(K_{\phi(p)}\cup K_{\phi(q)})\triangledown K_{\phi(n)+1}= P_{3}[K_{\phi(p)},K_{\phi(pq)+1},K_{\phi(q)} ] . \]
By Theorem \ref{mul of alpha n/n-1 of Zn}, $ \dfrac{n}{n-1} $ is the normalized Laplacian eigenvalue with multiplicity $ \phi(n) $. Again, by Theorems \ref{joined union} and \ref{zn}, we see that $ \dfrac{1}{q\phi(p)} $ and $ \dfrac{1}{q\phi(p)} $ are the normalized Laplacian eigenvalues of $ \mathcal{P}(\mathbb{Z}_n) $ with multiplicity $ \phi(p)-1 $ and $ \phi(q)-1 $ respectively. The remaining three normalized Laplacian eigenvalues are given by the following matrix
\begin{equation*} \label{Qmat pq}
\begin{pmatrix}
\dfrac{\phi(n)+1}{q\phi(p)} & \dfrac{-(\phi(n)+1)}{\sqrt{q\phi(p)(q\phi(p)+\phi(q))}} & 0\\
\dfrac{-\phi(p)}{\sqrt{q\phi(p)(q\phi(p)+\phi(q))}} &\dfrac{\phi(p)+\phi(q)}{q\phi(p)+\phi(q)} & \dfrac{-\phi(q)}{\sqrt{p\phi(q)(q\phi(p)+\phi(q))}}\\
0 & \dfrac{-(\phi(n)+1)}{\sqrt{p\phi(q)(q\phi(p)+\phi(q))}} & \dfrac{\phi(n)+1}{p\phi(q)}
\end{pmatrix}.
\end{equation*}$\hfill\Box$

By Corollaries \ref{NL spectra of p^z} and \ref{NL spectra of pq}, we have the following proposition.
\begin{proposition}
Equality holds in Theorem \eqref{mul of alpha n/n-1 of Zn}, if n is some prime or product of two primes.
\end{proposition}

The following theorem \cite{mehreen} shows that the power graph of a finite cyclic group $ \mathbb{Z}_{n} $ can be written as the joined union each of whose components are cliques.

\begin{theorem}\label{mehreen}
If $ \mathbb{Z}_{n} $ is a finite cyclic group of order $n\ge 3$, then the power graph $\mathcal{P}(\mathbb{Z}_{n})$ is given by
\begin{equation*}
\mathcal{P}(\mathbb{Z}_{n}) =K_{\phi(n)+1}\triangledown \mathbb{G}_{n}[K_{\phi(d_{1})},K_{\phi(d_{2})},\dots,K_{\phi(d_{t})}],
\end{equation*}
where $\mathbb{G}_{n}$ is the graph of order $t$ defined above.
\end{theorem}

Using Theorem \ref{joined union} and its consequences, we can compute the normalized Laplacian eigenvalues of $ \mathcal{P}(\mathbb{Z}_{n}) $ in terms of the adjacency spectrum of $ K_{\omega} $ and zeros of the characteristic polynomial of the auxiliary matrix.

We form a connected graph $ H=K_{1}\triangledown \mathbb{G}_{n} $ which is of diameter at most two if $ \mathbb{G}_{n} $ is not complete, otherwise its diameter is $ 1 $. In the following result, we compute the normalized Laplacian eigenvalues of the power graph of $ \mathbb{Z}_{n} $ by using Theorems \ref{joined union} and \ref{mehreen}.

\begin{theorem}\label{zn}
The normalized Laplacian eigenvalues of $ \mathcal{P}(\mathbb{Z}_n) $ are
\begin{equation*}
\left\lbrace \left (\dfrac{n}{n-1}\right )^{(\phi(n))}, \left (\dfrac{\phi(d_{1})+\alpha_{2}}{d_{1}+\alpha_{2}-1}\right )^{[\phi(d_{1})-1]},\dots, \left (\dfrac{\phi(d_{t})+\alpha_{r+1}}{d_{t}+\alpha_{t+1}-1}\right )^{[\phi(d_{t})-1]} \right\rbrace
\end{equation*}
and the $ t+1 $ eigenvalues of the following matrix $ \mathcal{M} $
\begin{equation}\label{quotient matrix of Z_n}
\mathcal{M} =\begin{pmatrix}
\dfrac{n-1-\phi(n)}{n-1}& \dfrac{-\phi(d_{1})a_{12}}{\sqrt{(\phi(n)+\alpha_{1})(r_{2}+\alpha_{2})}}&\ldots & \dfrac{-\phi(d_{t})a_{1(t+1)}}{\sqrt{(\phi(n)+\alpha_{1})(r_{t+1}+\alpha_{t+1})}}\\
\dfrac{-\phi(d_{1})a_{21}}{\sqrt{(r_{2}+\alpha_{2})(\phi(n)+\alpha_{1})}}&\dfrac{\alpha_{2}}{\alpha_2+r_2}& \ldots & \dfrac{-\phi(d_{t})a_{2(t+1)}}{\sqrt{(r_{2}+\alpha_{2})(r_{t+1}+\alpha_{t+1})}}\\
\vdots &\vdots  &\ddots &\vdots\\
\dfrac{-\phi(d_{1})a_{(t+1)1}}{\sqrt{(r_{n}+\alpha_{t+1})(\phi(n)+\alpha_{1})}}& \dfrac{-\phi(d_{2})a_{(t+1)2}}{\sqrt{(r_{n}+\alpha_{t+1})(r_{2}+\alpha_{2})}}&\ldots &\dfrac{\alpha_{t+1}}{\alpha_{t+1}+r_{t+1}}
\end{pmatrix},
\end{equation}
where, for $i\ne j$,
\begin{equation*}
a_{ij}=\begin{cases}
1, & v_i\sim v_j\\
0, & v_i\nsim v_j

\end{cases}
\end{equation*}
and $ r_{i}=\phi(d_{i})-1, ~ \text{for} ~i=2,\dots, t+1. $
\end{theorem}
\textbf{Proof.} Let $ \mathbb{Z}_{n} $ be a finite cyclic group of order $ n $. Since the identity element $ 0 $ and the $ \phi(n) $ generators of the group $ \mathbb{Z}_{n} $ are adjacent to every other vertex of $ \mathcal{P}(\mathbb{Z}_{n}) $, therefore, by Theorem \ref{mehreen}, we have
\begin{equation*}
\mathcal{P}(\mathbb{Z}_{n}) =K_{\phi(n)+1}\triangledown \mathbb{G}_{n}[K_{\phi(d_{1})},K_{\phi(d_{2})},\dots,K_{\phi(d_{t})}]= H[K_{\phi(n)+1},K_{\phi(d_{1})},K_{\phi(d_{2})},\dots,K_{\phi(d_{t})}],
\end{equation*}
where $ H=K_{1}\triangledown \mathbb{G}_{n} $ is the graph with vertex set $ \{v_{1},\dots, v_{t+1}\} $. Taking $ G_{1}=K_{\phi(n)+1}$ and $ G_{i}=K_{\phi(d_{i-1})} $, for $ i=2,\dots,t+1 $, in  Theorem \ref{joined union} and using the fact that the adjacency spectrum of $K_{\omega}$ consists of the eigenvalue $\omega-1$ with multiplicity $1$ and the eigenvalue $-1$ with multiplicity $\omega-1$, it follows that
\begin{equation*}
1-\dfrac{1}{r_{1}+\alpha_{1}}\lambda_{1k}A(G_{1}) =1-\dfrac{1}{r_{1}+\alpha_{1}}(-1)= 1+\dfrac{1}{\phi(n)+n-\phi(n)-1}=\dfrac{n}{n-1}
\end{equation*}
is a normalized eigenvalue of $ \mathcal{P}(\mathbb{Z}_n) $ with multiplicity $ \phi(n)$. Note that we have used the fact that vertex $ v_{1} $ of graph $ H $ is adjacent to every other vertex of $ H $ and $ \alpha_1=\sum\limits_{ d|n, d\ne 1,n}\phi(d)=n-1-\phi(n)$, as $ \sum\limits_{ d|s}\phi(d) =s.$ Similarly, we can show that $ \dfrac{\phi(d_{1})+\alpha_{2}}{\phi(d_{1})+\alpha_{2}-1},\dots, \dfrac{\phi(d_{t})+\alpha_{t+1}}{\phi(d_{t})+\alpha_{t+1}-1} $ are the normalized Laplacian eigenvalues of $ \mathcal{P}(\mathbb{Z}_n) $ with multiplicities $ \phi(d_{1})-1, \dots, \phi(d_{t})-1 $, respectively.  The remaining normalized Laplacian eigenvalues are the eigenvalues of the quotient matrix $ \mathcal{M} $ given by \eqref{quotient matrix of Z_n}. $\hfill\Box$

From Theorem \ref{zn}, it is clear that all the normalized Laplacian eigenvalues of the power graph $ \mathcal{P}(\mathbb{Z}_{n})$ are completely determined except the $t+1$ eigenvalues, which are the eigenvalues of the matrix $\mathcal{M}$ in Equation \eqref{quotient matrix of Z_n}. Further, it is also clear that the matrix $\mathcal{M}$  depends upon the structure of the graph $\mathbb{G}_{n} $, which is not known in general. However, if we give some particular value to $n$, then it may be possible to know the structure of graph $\mathbb{G}_{n}$ and hence about the matrix $\mathcal{M}$. This information may be helpful to determine the $t+1$ remaining normalized Laplacian eigenvalues of the power graph $ \mathcal{P}(\mathbb{Z}_{n})$.

We discuss some particular cases of Theorem \ref{zn}.\\

\begin{figure}
\centering

\begin{tikzpicture}
\draw[fill=black] (0,0) circle (2pt);\draw[fill=black] (2,0) circle (2pt);\draw[fill=black] (4,0) circle (2pt);
\draw[fill=black] (-1,-1) circle (2pt);\draw[fill=black] (1,-1) circle (2pt);\draw[fill=black] (3,-1) circle (2pt);

\draw[thin] (0,0)--(-1,-1); \draw[thin] (0,0)--(1,-1);\draw[thin] (2,0)--(-1,-1);
\draw[thin] (2,0)--(3,-1); \draw[thin] (4,0)--(1,-1);\draw[thin] (4,0)--(3,-1);

\node at (0,.3) {$ p $};\node at (2,.3) {$ q $};\node at (4,.3) {$ r $};
\node at (-1,-1.3) {$ pq $};\node at (1,-1.3) {$ pr $};\node at (3,-1.3) {$ qr $};
\end{tikzpicture}\quad
\begin{tikzpicture}
\draw[fill=black] (0,0) circle (2pt);\draw[fill=black] (2,0) circle (2pt);\draw[fill=black] (4,0) circle (2pt);
\draw[fill=black] (-1,-1) circle (2pt);\draw[fill=black] (1,-1) circle (2pt);\draw[fill=black] (3,-1) circle (2pt);
\draw[fill=black] (2.9,-2.3) circle (2pt);

\draw[thin] (0,0)--(2.9,-2.3);\draw[thin] (2,0)--(2.9,-2.3);\draw[thin] (3,-1)--(2.9,-2.3);
\draw[thin] (2,0)--(2.9,-2.3); \draw[thin] (4,0)--(2.9,-2.3);\draw[thin] (-1,-1)--(2.9,-2.3);

\draw[thin] (0,0)--(-1,-1); \draw[thin] (0,0)--(1,-1);\draw[thin] (2,0)--(-1,-1);
\draw[thin] (2,0)--(3,-1); \draw[thin] (4,0)--(1,-1);\draw[thin] (4,0)--(3,-1);

\node at (0,.3) {$ p $};\node at (2,.3) {$ q $};\node at (4,.3) {$ r $};
\node at (-1,-1.3) {$ pq $};\node at (1,-1.3) {$ pr $};\node at (3,-0.7) {$ qr $}; \node at (2.9,-2.7) {$ K_{1} $};
\end{tikzpicture}

\caption{Divisor graph $ \mathbb{G}_{pqr} $ and $ H=K_{1}\triangledown \mathbb{G}_{pqr} $. }
\end{figure}

Now, let $ n=pqr $, where $p,~ q,~ r$  with $ p<q<r$ are primes. From the definition of $\mathbb{G}_{n}$, the  vertex set and edge set of $\mathbb{G}_{n}$ are $\{ p, q, r, pq, pr, qr \}$ and $\{(p,pq),(p,pr),(q,pq),(q,qr),(r,pr),(r,qr)\} $ respectively, and is shown in Figure $ (1) $. Let $ H=K_{1}\triangledown \mathbb{G}_{n}.$ Then
\begin{equation*}
 \mathcal{P}(\mathbb{Z}_n)= H[K_{\phi(n)+1}, K_{\phi(p)},K_{\phi(q)},K_{\phi(r)},K_{\phi(pq)},K_{\phi(pr)}, K_{\phi(qr)}].
\end{equation*}
By Theorem \ref{joined union}, we have \begin{align*}
\big (\alpha_{1},&\alpha_{2}, \alpha_{3},\alpha_{4},\alpha_{5},\alpha_{6},\alpha_{7}\big)=\big(n-\phi(n)-1,\phi(n)+1+\phi(pq)+\phi(pr),\phi(n)+1+\phi(pq)+\phi(qr)),\\
\phi(n)&+1+\phi(pr)+\phi(qr),\phi(n)+1+\phi(p)+\phi(q),\phi(n)+1+\phi(p)+\phi(r),\phi(n)+1+\phi(q)+\phi(r)\big ) \\
\text{and}&\quad \big (\alpha_{1}+r_{1},\alpha_{2}+r_{2}, \alpha_{3}+r_{3},\alpha_{4}+r_{4},\alpha_{5}+r_{5},\alpha_{6}+r_{6},\alpha_{7}+r_{7}\big)=\big(n-1,\phi(n)+\phi(p)+\phi(pq)\\
+\phi&(pr),\phi(n)+\phi(q)+\phi(pq)+\phi(qr)),\phi(n)+\phi(r)+\phi(pr)+\phi(qr),\phi(n)+\phi(pq)+\phi(p)+\phi(q),\\
\phi&(n)+\phi(pr)+\phi(p)+\phi(r),\phi(n)+\phi(qr)+\phi(q)+\phi(r)\big ).
\end{align*}
Also, by Theorem \ref{mul of alpha n/n-1 of Zn}, $ \dfrac{n}{n-1} $ is the normalized Laplacian eigenvalue with multiplicity $ \phi(n) $.  Using the above information and Theorem \ref{zn}, the second distinct normalized Laplacian eigenvalue is $ 1+\dfrac{1}{r_{2}+\alpha_{2}}=1+\dfrac{1}{\phi(n)+\phi(p)+\phi(pq)+\phi(pr)} $ with multiplicity $ \phi(p)-1 $. In a similar way, we see that the other eigenvalues are
\begin{align*}
&1+\dfrac{1}{\phi(n)+\phi(q)+\phi(pq)+\phi(qr)},1+\dfrac{1}{\phi(n)+\phi(r)+\phi(pr)+\phi(qr)},1+\dfrac{1}{\phi(n)+\phi(pq)+\phi(q)+\phi(q)},\\
&1+\dfrac{1}{\phi(n)+\phi(pr)+\phi(p)+\phi(r)},~1+\dfrac{1}{\phi(n)+\phi(qr)+\phi(q)+\phi(r)}
\end{align*}  with multiplicities $ \phi(q)-1,~ \phi(r)-1,~ \phi(pq)-1,~ \phi(pr)-1,~ \phi(qr)-1$, respectively. The remaining $7$ eigenvalues are given by the following matrix
\begin{equation*}\label{pqr}
\begin{pmatrix}
z_{1}& -\phi(p)c_{12} & -\phi(q)c_{13} & -\phi(r)c_{14}& -\phi(pq)c_{15}& -\phi(pr)c_{16}& -\phi(qr)c_{17}\\
(\phi(n)+1)c_{21} & z_{2} & 0 & 0& -\phi(pq)c_{25} & -\phi(pr)c_{26}& 0 \\
(\phi(n)+1)c_{31} & 0 & z_{3} & 0& -\phi(pq)c_{35} & 0 & -\phi(qr)c_{37} \\
(\phi(n)+1)c_{41} & 0 & 0 & z_{4} & 0&  -\phi(pr)c_{46} & -\phi(qr)c_{47} \\
(\phi(n)+1)c_{51} & -\phi(p)c_{25} & -\phi(q)c_{35} & 0 & z_{5} & 0  & 0  \\
(\phi(n)+1)c_{61} & -\phi(p)c_{26} & 0 & -\phi(r)c_{64} & 0 & z_{6} & 0  \\
(\phi(n)+1)c_{71} & 0 & -\phi(q)c_{75} & -\phi(r)c_{74} & 0 & 0 &  z_{7}
\end{pmatrix},
\end{equation*}
where, \begin{align*}
z_{1}=&\frac{n-\phi(n)-1}{n-1}, ~ z_{2}= \frac{\phi(n)+1+\phi(pq)+\phi(pr)}{\phi(n)+\phi(p)+\phi(pq)+\phi(pr)},~ z_{3}= \frac{\phi(n)+1+\phi(pq)+\phi(qr)}{\phi(n)+\phi(q)+\phi(pq)+\phi(qr)},\\
&z_{4}=\dfrac{\phi(n)+1+\phi(pr)+\phi(qr)}{\phi(n)+\phi(r)+\phi(pr)+\phi(qr)}, ~ z_{5}= \dfrac{\phi(n)+1+\phi(p)+\phi(q)}{\phi(n)+\phi(pq)+\phi(p)+\phi(q)},\\
& z_{6}= \dfrac{\phi(n)+1+\phi(p)+\phi(r)}{\phi(n)+\phi(pr)+\phi(p)+\phi(r)}, z_{7}=\dfrac{\phi(n)+1+\phi(q)+\phi(r)}{\phi(n)+\phi(qr)+\phi(q)+\phi(r)}\\
\text{and}~c_{ij}&=c_{ji}=\dfrac{1}{\sqrt{(r_{i}+\alpha_{i})(r_{j}+\alpha_{j})}}.
\end{align*}\\

Next, we discuss the normalized Laplacian eigenvalues of the finite cyclic group $\mathbb{Z}_{n}$, with $n=p^{n_{1}}q^{n_{2}},$ where $ p<q $ are primes and $ n_{1}\leq  n_{2}$ are positive integers. We consider the case when both $ n_{1} $ and $ n_{2} $ are even, and the case when they are odd can be discussed similarly.

\begin{theorem}
Let $ \mathcal{P}(\mathbb{Z}_{p^{n_{1}}q^{n_{2}}}) $ be the power graph of the finite cyclic group $\mathbb{Z}_{p^{n_{1}}q^{n_{2}}}$ of order $ n={p^{n_{1}}q^{n_{2}}} $, where $p<q$ are primes and $n_{1}=2m_{1}\leq n_{2}=2m_{2}$ are even positive integers. Then the normalized Laplacian eigenvalues of $ \mathcal{P}(\mathbb{Z}_{p^{n_{1}}q^{n_{2}}}) $ consists of the eigenvalues
\begin{align*}
&\left (\dfrac{n}{n-1}\right )^{[\phi(n)]}, \left (\dfrac{n-q^{n_{2}}+1}{n-q^{n_{2}}}\right )^{[\phi(p)-1]},\\
&\qquad \vdots\\
&\left (\dfrac{p^{m_{1}-1}+q^{n_{2}}(p^{n_{1}}-p^{m_{1}-1})}{p^{m_{1}-1}+q^{n_{2}}(p^{n_{1}}-p^{m_{1}-1})}-1\right )^{[\phi(p^{m_{1}})-1]},\\
&\qquad\vdots\\
&\left (\dfrac{p^{n_{1}-1}+q^{n_{2}}\phi(p^{n_{1}})}{p^{n_{1}-1}+q^{n_{2}}\phi(p^{n_{1}})-1}\right )^{\left [\phi\left (p^{n_{1}}\right )-1\right ]},\left (\dfrac{n-p^{n_{1}}+1}{n-p^{n_{1}}}\right )^{\left [\phi\left (q\right )-1\right ]},\\
&\qquad\vdots\\
&\left (\dfrac{q^{m_{2}-1}+p^{n_{1}}(q^{n_{2}}-q^{m_{2}-1})}{q^{m_{2}-1}+p^{n_{1}}(q^{n_{2}}-q^{m_{2}-1})-1}\right )^{\left [\phi\left (q^{m_{2}}\right )-1\right ]},\\
&\qquad\vdots\\
&\left (\dfrac{q^{n_{2}-1}+p^{n_{1}}\phi(q^{n_{2}})}{q^{n_{2}-1}+p^{n_{1}}\phi(q^{n_{2}})-1}\right )^{\left [\phi\left (q^{n_{2}}\right )-1\right ]},\left (\dfrac{\phi(p)+\phi(q)+(q^{n_{2}}-1)(p^{n_{1}}-1)+1}{\phi(p)+\phi(q)+(q^{n_{2}}-1)(p^{n_{1}}-1)}\right )^{\left [\phi\left (pq\right )-1\right ]},\\
&\qquad\vdots\\
&\left (\dfrac{q^{n_{2}}(p^{n_{1}}-1)+q^{m_{2}}-q^{m_{2}-1}(p^{n_{1}}-p)}{q^{n_{2}}(p^{n_{1}}-1)+q^{m_{2}}-q^{m_{2}-1}(p^{n_{1}}-p)-1}\right )^{\left [\phi\left (pq^{m_{2}}\right )-1\right ]},\\
&\qquad\vdots\\
& \left (\dfrac{pq^{n_{2}}+\phi(p^{n_{1}})(q^{n_{2}}-q)}{pq^{n_{2}}+\phi(p^{n_{1}})(q^{n_{2}}-q)-1} \right )^{\left [\phi\left (pq^{n_{2}}\right )-1\right ]},\\
&\qquad\vdots\\
&\left (\dfrac{p^{m_{1}}+p^{n_{1}}(q^{n_{2}}-1)-p^{m_{1}-1}(q^{n_{2}}-q)}{p^{m_{1}}+p^{n_{1}}(q^{n_{2}}-1)-p^{m_{1}-1}(q^{n_{2}}-q)-1} \right )^{\left [\phi\left (p^{m_{1}}q\right )-1\right ]},\\
&\qquad\vdots\\
&\left (\dfrac{n+p^{m_{1}}q^{m_{2}}+p^{m_{1}-1}q^{m_{2}-1}-\phi(p^{m_{1}}q^{m_{2}})-p^{n_{1}}q^{m_{2}-1}-p^{m_{1}-1}q^{n_{2}}}{n+p^{m_{1}}q^{m_{2}}+p^{m_{1}-1}q^{m_{2}-1}-\phi(p^{m_{1}}q^{m_{2}})-p^{n_{1}}q^{m_{2}-1}-p^{m_{1}-1}q^{n_{2}}-1} \right )^{\left [\phi\left (p^{m_{1}}q^{m_{2}}\right )-1\right ]},\\
&\qquad\vdots\\
&\left (\dfrac{p^{m_{1}}q^{n_{2}}+\phi(q^{n_{2}})(p^{n_{1}}-p^{m_{1}})}{p^{m_{1}}q^{n_{2}}+\phi(q^{n_{2}})(p^{n_{1}}-p^{m_{1}})-1} \right )^{\left [\phi\left (p^{m_{1}}q^{n_{2}}\right )-1\right ]},\\
&\qquad\vdots\\
&\left (\dfrac{p^{n_{1}}q+\phi(p^{n_{1}})(q^{n_{2}}-q)}{p^{n_{1}}q+\phi(p^{n_{1}})(q^{n_{2}}-q)-1} \right )^{\left [\phi\left (p^{n_{1}}q\right )-1\right ]},\\
&\qquad\vdots\\
& \left ( \dfrac{p^{n_{2}}q^{m_{1}}+\phi(p^{n_{1}})(q^{n_{2}}-q^{m_{1}})}{p^{n_{2}}q^{m_{1}}+\phi(p^{n_{1}})(q^{n_{2}}-q^{m_{1}})-1} \right )^{\left [\phi\left (p^{n_{1}}q^{m_{2}}\right )-1\right ]},\\
&\qquad\vdots\\
&\left ( \dfrac{p^{n_{1}}q^{n_{2}-1}+\phi(n)}{p^{n_{1}}q^{n_{2}-1}+\phi(n)-1} \right )^{\left [\phi\left (p^{n_{1}}q^{n_{2}-1}\right )-1\right ]}
\end{align*}
and the remaining eigenvalues are given by matrix \eqref{quotient matrix of Z_n}.
\end{theorem}
\textbf{Proof.} Suppose that $n=p^{n_{1}}q^{n_{2}}$, where $ n_{1}=2m_{1} $ and $ n_{2}=2m_{2} $ are even with $ n_{1}\leq n_{2} $ and $ m_{1} $ and $ m_{2} $ are positive integers. Since the total number of divisors of $ n $ are $ (n_{1}+1)(n_{2}+1) $, so the order of $ \mathbb{G}_{p^{n_{1}}q^{n_{2}}} $ is $ (n_{1}+1)(n_{2}+1)-2 $. The proper divisor set of $n$ is
\begin{align*}
D(n)=\bigg\{p,p^{2},\cdots,p^{m_{1}}&,\dots,p^{n_{1}}, q,q^2,\dots,q^{m_{2}},\dots,q^{n_{2}}, pq, pq^2,\dots,pq^{m_{2}},\dots,pq^{n_{2}}, \cdots,p^{m_{1}}q, p^{m_{1}}q^2,\\
\dots, &p^{m_{1}}q^{m_{2}},\dots,p^{m_{1}}q^{n_{2}},\dots,p^{n_{1}}q, p^{n_{1}}q^2,\dots,p^{n_{1}}q^{m_{2}},\dots,p^{n_{1}}q^{n_{2}-1}\bigg\}.
\end{align*}
By the definition of graph $ \mathbb{G}_{n} $, we see that $ p$ is not adjacent to $ p, q,q^{2},\cdots,q^{m_{2}},\cdots,q^{n_{2}} $. So we write adjacency of vertices in terms of iterations and avoid divisors outside the  set $ D(n) $. Thus, we observe that
\begin{align*}
p&\sim p^{i},p^{j}q^{k},~ \text{for}~ i=2,3,\dots,n_{1},~  j=1,2,\dots,n_{1},~ k=1,2,\dots,n_{2},  \\
&\vdots\\
p^{m_{1}}&\sim p^{i},p^{j}q^{k},~ \text{for}~ i=1,2,\dots,n_{1}, ~i\neq m_{1},~ j=m_{1},\dots,n_{1},~ k=1,2,\dots,n_{2},\\
&\vdots\\
p^{n_{1}}&\sim p^{i},p^{n_{1}}q^{j},~ \text{for}~ i=1,2,\dots,n_{1}-1,~  j=1,2,\dots,n_{2}-1,\\ 
q&\sim q^{i},p^{j}q^{k},~ \text{for}~ i=2,3,\dots,n_{2}, ~ j=1,2,\dots,n_{1},~ k=1,2,\dots,n_{2},\\
& \vdots\\
q^{m_{2}}&\sim q^{i},p^{j}q^{k},~ \text{for}~ i=1,2,3,\dots,n_{1}, ~i\neq m_{2}, ~j=1,2,3,\dots,n_{1}, ,~ k=m_{2},\dots,n_{2}\\
&\vdots\\
q^{n_{2}}&\sim q^{i},p^{j}q^{n_{2}},~ \text{for}~ i=1,2,3,\dots,n_{2}-1,~ j=1,2,3,\dots,n_{1}-1,\\ 
pq&\sim p, q, p^{i}q^{j},~ \text{for } ~ i=1,2,3,\dots,n_{1},~ j=1,2,3,\dots,n_{2},\\
& \vdots\\
pq^{m_{2}}&\sim  p, q^{i},pq^{j},p^{k}q^{k} ~ \text{for}~ i=1,2,3,\dots,m_{2}, ~ j=1,2,3,\dots,n_{2},~ j\neq m_{2},~ k=2,3,\dots,n_{1},\\
&l=m_{2},\dots,n_{2},\\
&\vdots\\
pq^{n_{2}}&\sim  p, q^{i}, pq^{j},p^{k}q^{n_{2}}, ~\text{for}~ i=1,2,\dots,n_{2},~ j=1,2,\dots,n_{2}-1,~k=2,3,\dots,n_{1}-1,\\ 
& \vdots\\
p^{m_{1}}q &\sim  p^{i}, q,p^{m_{1}}q^{j},p^{k}q,p^{l}q^{m} ~\text{for}~ i=1,2,3,\dots,m_{1}, ~ j=2,3,\dots,n_{2},~k=1,2,\dots,m_{1}-1,\\
&l=m_{1}+1,\dots,n_{1},~ m=1,2,\dots,n_{2},\\
& \vdots\\
p^{m_{1}}q^{m_{2}} &\sim  p^{i}, q^{j}, p^{k}q^{l}~ \text{for}~ i=1,2,\dots,m_{1}, ~j=1,2,\dots,m_{2}, ~k=1,2,\dots,n_{1}, ~l=1,2,\dots,n_{2},\\
&\vdots\\
p^{m_{1}}q^{n_{2}} &\sim  p^{i}, q^{j},p^{k}q^{n_{2}}, p^{i}q^{j}~ \text{for}~ i=1,2,\dots,m_{1},~j=1,2,\dots,n_{2},~k=m_{1}+1,m_{1}+2,\dots,n_{1}-1,\\
&\vdots\\
p^{n_{1}}q &\sim  p^{i}, q ,p^{j}q, p^{n}q^{k}~ \text{for}~ i=1,2,\dots,n_{1},~j=1,2,\dots,n_{1}-1,~k=2,3\dots,n_{2}-1,\\
&\vdots\\
p^{n_{1}}q^{m_{2}} &\sim  p^{i}, q^{j} ,p^{n_{1}}q^{k}, p^{i}q^{j}~ \text{for}~ i=1,2,\dots,n_{1},~j=1,2,\dots,m_{2},~k=m_{2}+1,m_{2}+2\dots,n_{2}-1,\\
&\vdots\\
p^{n_{1}}q^{n_{2}-1} &\sim  p^{i}, q^{j} , p^{i}q^{j}~ \text{for}~ i=1,2,\dots,n_{1},~j=1,2,\dots,n_{2}-1.\\
\end{align*}
Therefore, by Theorem \ref{mehreen}, we have \begin{align*}
\mathcal{P}(\mathbb{Z}_{n})= &K_{\phi(n)+1}\triangledown \mathbb{G}_{n}[K_{\phi(p)},\dots,K_{\phi(p^{m_{1}})},\dots,K_{\phi(p^{n_{1}})},K_{\phi(q)},\dots,K_{\phi(q^{m_{2}})},\dots,K_{\phi(q^{n_{2}})},K_{\phi(pq)},\dots,\\
&K_{\phi(pq^{m_{2}})},\dots,K_{\phi(pq^{n_{2}})},\dots,K_{\phi(p^{m_{1}}q)},\dots,K_{\phi(p^{m_{1}}q^{m_{2}})},\dots,K_{\phi(p^{m_{1}}q^{n_{2}})},\dots,K_{\phi(p^{n_{1}}q)},\dots,\\
&K_{\phi(p^{n_{1}}q^{m_{2}})},\dots,K_{\phi(p^{n_{1}}q^{n_{2}-1})}].\\
\end{align*}
Now, by using Theorem \ref{joined union}, we calculate the values of $ \alpha_{i} $'s and $ r_{i}+\alpha_{i}=r_{i}^{'} $'s. We recall some number theory identities, like $\phi(xy)=\phi(x)\phi(y),$ provided that $(x,y)=1$, $\sum\limits_{i=1}^{k}\phi(p^i)=p^k-1$ and $ \sum\limits_{d|s}\phi(d)=s $. Using this information and definition of $ \alpha_{i} $'s, we have
\begin{equation*}
\alpha_{1}=\sum\limits_{1,n\neq d|n}\phi(d)=n-1-\phi(n)
\end{equation*}
and
\begin{equation*}
\begin{split}
\alpha_{2}&= \phi(p^{2})+\dots+\phi(p^{m_{1}})+\dots+\phi(p^{n_{1}})+\phi(pq)+\dots+\phi(pq^{m_{2}})+\dots\phi(pq^{n_{2}})\\&
+\phi(p^{m_{1}}q)+\dots+\phi(p^{m_{1}}q^{m_{2}})+\dots\phi(p^{m_{1}}q^{n_{2}})+\dots+\phi(p^{n_{1}}q)
+\dots+\phi(p^{n_{1}}q^{m_{2}})\\&
+\dots\phi(p^{n_{1}}q^{n_{2}-1})+\phi(n)+1\\&
=\sum\limits_{1,p,n\neq d|n}\phi(d)-[\phi(q)+\dots+\phi(q^{m_{1}})+\dots+\phi(q^{n_{1}})]+\phi(n)+1\\&
=n-1-\phi(p)-\phi(n)-[q^{n_{2}}-1]+\phi(n)+1=n-\phi(p)-q^{n_{2}}+1.
\end{split}
\end{equation*}
Proceeding in the same way as above, other $ \alpha_{i} $'s are
\begin{align*}
\alpha_{3}=&q^{n_{2}}(p^{n_{1}}-p)+p-\phi(p^{2}),\\
& \vdots\\
\alpha_{m_{1}+1}=& p^{m_{1}-1}+q^{n_{2}}(p^{n_{1}}-p^{m_{1}-1}) -\phi(p^{m_{1}}),\\
\vdots\\
\alpha_{n_{1}+1}=&p^{n_{1}-1}+\phi(p^{n_{1}})(q^{n_{2}-1}-1),~\alpha_{n_{1}+2} =n-\phi(q)-p^{n_{1}}+1,\\
&\vdots\\
\alpha_{n_{1}+m_{1}+1}=&q^{m_{2}-1}+p^{n_{1}}(q^{n_{2}}-q^{m_{2}-1})-\phi(q^{m_{2}}),\\
\vdots\\
\alpha_{n_{1}+n_{2}+1} =&q^{n_{2}-1}+\phi(q^{n_{2}})(p^{n_{1}}-1),\\
\alpha_{n_{1}+n_{2}+2}=&\phi(p)+\phi(q)+1-\phi(pq)+(q^{n_{2}}-1)(p^{n_{1}}-1),\\
&\vdots\\
\alpha_{n_{1}+n_{2}+m_{1}+1}=&q^{n_{2}}(p^{n_{1}}-1)+q^{m_{2}}-q^{m_{2}-1}(p^{n_{1}}-p)-\phi(pq^{m_{2}}),\\
&\vdots\\
\alpha_{2n_{1}+n_{2}+1}=&pq^{n_{2}}-\phi(pq^{n_{2}})+\phi(p^{n_{1}})(q^{n_{2}}-q),\\
&\vdots\\
\alpha_{m_{1}n_{2}+n_{1}+2}=&p^{m_{1}}-\phi(p^{m_{1}}q)+p^{n_{1}}(q^{n_{2}}-1) -p^{m_{1}-1} (q^{v_{1}}-q),\\
&\vdots\\
 \alpha_{m_{1}n_{2}+n_{1}+m_{1}+1}=&n+p^{m_{1}-1}q^{m_{2}-1}+p^{m_{1}}q^{m_{2}}-2\phi(p^{m_{1}}q^{m_{2}})-q^{n_{2}}p^{m_{1}-1}-p^{n_{2}}q^{m_{2}-1},\\
&\vdots\\
\alpha_{(m_{1}+1)n_{2}+n_{1}+1}=&p^{m_{1}}q^{n_{2}}-\phi(p^{m_{1}}q^{n_{2}})+\phi(q^{n_{2}})(p^{n_{1}}-p^{m_{1}}),\\
&\vdots,\\
\alpha_{n_{1}n_{2}+n_{1}+2}=&p^{n_{1}}q+\phi(p^{n_{1}})(q^{n_{2}}-q)-\phi(p^{n_{1}}q),\\
&\vdots,\\
\alpha_{n_{1}n_{2}+n_{1}+m_{1}+1}=&p^{n_{1}}q^{m_{2}}+\phi(p^{n_{1}})(q^{n_{2}}-q^{m_{2}})-\phi(p^{n_{1}}q^{m_{2}}),\\
&\vdots\\
\alpha_{(n_{1}+1)n_{2}+n_{1}}=&p^{n_{1}}q^{n_{2}-1}+\phi(n)-\phi(p^{n_{1}}q^{n_{2}-1}).
\end{align*}
Also, value of $r_{i}+\alpha_{i}= r_{i}^{'} $'s are given by
\begin{align*}
r_{1}^{'}=&n-1,~ r_{2}^{'}=n-q^{n_{2}},\\
&\vdots\\
r_{m_{1}+1}^{'}=&p^{m_{1}-1}+q^{n_{2}}(p^{n_{1}}-p^{m_{1}-1})-1,\\
&\vdots\\
r_{n_{1}+1}^{'}=&p^{n_{1}-1}+\phi(p^{n_{1}})q^{n_{2}}-1, ~r_{n_{1}+2}^{'}=n-p^{n_{1}},\\
&\vdots\\
r_{n_{1}+m_{1}+1}^{'}=&q^{m_{2}-1}+p^{n_{1}}(q^{n_{2}}-q^{m_{2}-1})-1,\\
&\vdots\\
r_{n_{1}+n_{2}+1}^{'}=&q^{n_{2}-1}+\phi(q^{n_{2}})p^{n_{1}}-1,~r_{n_{1}+n_{2}+2}^{'}=\phi(p)+\phi(q)+(q^{n_{2}}-1)(p^{n_{1}}-1),\\
&\vdots\\
r_{n_{1}+n_{2}+m_{1}+1}^{'}=&q^{n_{2}}(p^{n_{1}}-1)+q^{m_{2}}-q^{m_{2}-1}(p^{n_{1}}-p)-1,\\
&\vdots\\
r_{2n_{1}+n_{2}+1}^{'}=&pq^{n_{2}}+\phi(p^{n_{1}})(q^{n_{2}}-q)-1,\\
&\vdots\\
r_{m_{1}n_{2}+n_{1}+2}^{'}=&p^{m_{1}}+p^{n_{1}}(q^{n_{2}}-1) -p^{m_{1}-1} (q^{n_{1}}-q)-1,\\
&\vdots\\
r_{m_{1}N_{2}+N_{1}+m_{1}+1}^{'} =&n+p^{m_{1}-1}q^{m_{2}-1}+p^{m_{1}}q^{m_{2}}-\phi(p^{m_{1}}q^{m_{2}})-q^{n_{2}}p^{m_{1}-1}-p^{n_{2}}q^{m_{2}-1}-1,\\
&\vdots\\
r_{(m_{1}+1)n_{2}+n_{1}+1}^{'}=&p^{m_{1}}q^{n_{2}}+\phi(q^{n_{2}})(p^{n_{1}}-p^{m_{1}})-1,\\
&\vdots\\
r_{n_{1}n_{2}+n_{1}+2}^{'}=&p^{n_{1}}q+\phi(p^{n_{1}})(q^{n_{2}}-q)-1,\\
&\vdots\\
r_{n_{1}n_{2}+n_{1}+m_{1}+1}^{'}
=&p^{n_{1}}q^{m_{2}}+\phi(p^{n_{1}})(q^{n_{2}}-q^{m_{2}})-1,\\
&\vdots\\
r_{n_{1}n_{2}+n_{1}+m_{1}+1}^{'} =&p^{n_{1}}q^{n_{2}-1}+\phi(n)-1.
\end{align*}
We note that each of $ G_{i}=K_{i} $ and by Theorems \ref{mul of alpha n/n-1 of Zn} and \ref{joined union}, we get the desired eigenvalues as in the statement. By substituting the values of $ \alpha_{i} $'s, $ r_{i}^{'} $'s and using the adjacency relations, the remaining normalized Laplacian eigenvalues are the eigenvalues of matrix \eqref{quotient matrix of Z_n}.  $\hfill\Box$ \\

\noindent{\bf Acknowledgements.}  The research of S. Pirzada is supported by the SERB-DST research project number MTR/2017/000084.


\begin{thebibliography}{0}
\bibitem{survey} { J. Abawajy, A. Kelarev and M. Chowdhury, Power graphs: A survey, \emph{Electronic J. Graph Theory Appl.}} {\textbf{1(2)}} (2013) 125--147.

\bibitem{banerjee} {S. Banerjee and A. Adhikari, Signless Laplacian spectra of power graphs of certain finite groups}, \emph{AKCE Int. J. Graphs Comb.} DOI:10.1016/j.akcej.2019.03.009 (2019).
\bibitem{BH} {A. E. Brouwer, W. H. Haemers \emph{ Spectra of Graphs}}, Springer New York 2010.



\bibitem{cameron1} {P. J. Cameron and S. Ghosh, The power graphs of a finite group, \emph{ Dicrete Math.}} {\textbf{ 311(13)}} (2011) 1220--1222.
\bibitem{cavers} M. Cavers, The normalized Laplacian matrix and general Randic index of graphs, Thesis, University of Regina 2010.
\bibitem{sriparna} {S. Chattopadhyay and P. Panigrahi, On Laplacian spectrum of power graphs of finite cyclic and dihedral groups, \emph{ Linear Multilinear Algebra}} {\textbf{ 63(7)}} (2015) 1345--1355.
\bibitem{chung} {F. R. K. Chung},\emph{ Spectral Graph Theory} American Mathematical Society, Providence 1997.
\bibitem{sen} {I. Chakrabarty,  M. Ghosh and M. K. Sen, Undirected power graph of semigroups, \emph{ Semigroup Forum}} {\textbf{78}} (2009) 410--426.
\bibitem{tamiza} {T. T. Chelvan and M. Sattanathan,  Power graphs of finite abelian groups, \emph{ Algebra Discrete Math.}} {\textbf{16(1)}} (2013) 33--41.

\bibitem{cardosa} D. M. Cardoso, M. A. De Freitas, E. N. Martins and M Robbiano, Spectra of graphs obtained by a generalization of the join of graph operations, \emph {Discrete Math} {\bf 313} (2013) 733--741.

\bibitem{cds} D. M. Cvetkovi\'{c}, M. Doob and H. Sachs, \emph {Spectra of graphs. Theory and Applications,} Pure and Applied Mathematics, 87. Academic Press, Inc. New York 1980.
\bibitem{db} {K. C. Das, A. D. G\"{u}ng\"{o}r and \c{S}. Bozkurt, On the normalized Laplacian eigenvalues of graphs, \emph{ Ars Combinatoria}} \textbf{118} (2015) 143--154.
\bibitem{db1} {K. C. Das, S. Sun and I. Gutman, Normalized Laplacian eigenvalues and Randi\'c energy of graphs, \emph{ MATCH Comm. Math. Comp. Chem.}} \textbf{77(1)} (2017) 45--59.



\bibitem{asma1} {A. Hamzeh, Signless and normalized Laplacian spectrums of the power graph and its supergraphs of certain finite groups}, \emph{J. Indonesian Math. Soc.} {\textbf{ 24(1)}} (2018) 61--69.
\bibitem{asma} {A. Hamzeh and A. R. Ashrafi, Spectrum and L-spectrum of the power graph and its main supergraph for certain finite groups} \emph{Filomat} {\textbf{31(16)}} (2017) 5323--5334.

\bibitem{kelarev} {A. V. Kelarev and S. J. Quinn, Directed graphs and combinatorial properties of semigroups, \emph{J. Algebra}} {\textbf{251}} (2002) 16--26.
\bibitem{kelarev9} {A. V. Kelarev and S. J. Quinn, Graph algebras and automata, \textbf{257}, \emph{Marcel Dekker}} New York, 2003.


\bibitem{mehreen} Z. Mehranian, A. Gholami and A. R. Ashrafi, A note on the power graph of a finite group,\emph{ Int. J. Group Theory }{\textbf{ 5}} (2016) 1--10.
\bibitem{mehreen1} Z. Mehranian, A. Gholami and A. R. Ashrafi, The spectra of power graphs of certain finite groups,\emph{ Linear Multilinear Algebra}  {\textbf{65(5)}} (2016) 1003--1010.


\bibitem{roman}{ W. K. Nicholson, Introduction to abstract algebra, Fourth edition}, John Wiley and Sons, New Jersey (2012).


\bibitem{panda} {R. P. Panda, Laplacian spectra of power graphs of certain finite groups, \emph{ Graphs Combinatorics}} DOI:10.1007/s00373-019-02070-x (2019).
\bibitem{sp} S. Pirzada, \emph {An Introduction to Graph Theory}, Universities Press, Orient BlackSwan Hyderabad (2012).

\bibitem{bilal}  B. A. Rather, S. Pirzada and Z. Goufei, On distance Laplacian spectra of power graphs of certain finite groups, preprint.

\bibitem{DS} D. Stevanovi\'{c},  Large sets of long distance equienergetic graphs, \emph {Ars Math. Contemp.} \textbf{2(1)} (2009) 35--40.
\bibitem{sun} {S. Sun and K. C. Das, Normalized Laplacian spectrum of complete multipartite graphs, \emph{ Discrete Applied Math.}} \textbf{284} (2020) 234--245.

\bibitem{wu} B. F. Wu, Y. Y. Lou and C. X. He, Signless Laplacian and normalized Laplacian on the H-join operation of
graphs, \textit{Discrete Math. Algorithm. Appl.} \textbf{06} (2014) [13 pages] DOI:http://dx.doi.org/10.1142/S1793830914500463.



\end{thebibliography}
\end{document}